\documentstyle[12pt]{article}
\newcommand{\sect}[1]{\section{#1}\setcounter{equation}{0}}

\font\mbn=msbm10 scaled \magstep1
\font\mbs=msbm7 scaled \magstep1
\font\mbss=msbm5 scaled \magstep1
\newfam\mbff
\textfont\mbff=\mbn
\scriptfont\mbff=\mbs
\scriptscriptfont\mbff=\mbss\def\mbf{\fam\mbff}
\def\Re{{\mbf R}}

\def\Co{{\mbf C}}

\def\Di{{\mbf D}}
\def\Bo{{\mbf B}}
\def\N{{\mbf N}}
\newtheorem{Th}{Theorem}[section]
\newtheorem{Lm}[Th]{Lemma}
\newtheorem{C}[Th]{Corollary}

\newtheorem{Proposition}[Th]{Proposition}

\author{Alexander Brudnyi\thanks{Research supported in part by NSERC.\newline
}\\
Department of Mathematics and Statistics\\
University of Calgary, Calgary\\
Canada}
\title{Extension of Matrices with Entries in $H^{\infty}$
on Coverings of Riemann Surfaces of
Finite Type}
\date{} 
\begin{document} 
\maketitle
\begin{abstract}
{In the present paper continuing our work started in [Br1]-[Br5]
we prove an extension theorem for matrices with entries in the algebra of bounded
holomorphic functions defined on an unbranched covering of a Caratheodory
hyperbolic Riemann surface of finite type.} 
\end{abstract}
\sect{\hspace*{-1em}. Introduction.}
\quad
Let $X$ be a complex manifold and let $H^{\infty}(X)$ be the Banach algebra 
of bounded holomorphic functions on $X$
equipped with the supremum norm. We assume that $X$ is Caratheodory
hyperbolic, that is, the functions in 
$H^{\infty}(X)$ separate the points of $X$. 
The maximal ideal space ${\cal M}={\cal M}(H^{\infty}(X))$ is the set of all
nonzero linear multiplicative functionals on $H^{\infty}(X)$. Since the norm
of each $\phi\in {\cal M}$ is $\leq 1$, ${\cal M}$
is a subset of the closed unit ball of the
dual space $(H^{\infty}(X))^{*}$. It is a
compact Hausdorff space in the Gelfand topology (i.e., in the weak $*$ 
topology induced by $(H^{\infty}(X))^{*}$). Also, there is a
continuous embedding $i:X\hookrightarrow {\cal M}$ taking $x\in X$ to the 
evaluation homomorphism $f\mapsto f(x)$, $f\in H^{\infty}(X)$. 
The complement to
the closure of $i(X)$ in ${\cal M}$ is called the {\em corona}.
The {\em corona problem} is: given $X$ to determine whether the corona is 
empty. For example, according to Carleson's celebrated Corona Theorem
[C] this is true for $X$ being the open unit disk in $\Co$. (This was
conjectured by Kakutani in 1941.) Also, there are
non-planar Riemann surfaces for which the corona is non-trivial (see, e.g.,
[JM], [G], [BD], [L] and references therein). This is due to a structure that in a sense makes the surface seem higher dimensional. So there is a hope that the restriction to the Riemann sphere might prevent this obstacle.
However, the general problem for planar
domains is still open, as is the problem in several variables for the ball
and polydisk. (In fact, there are no known examples of domains in $\Co^{n}$,
$n\geq 2$, without corona.) At present, the strongest
corona theorem for planar domains is due to Moore [M]. It states that the corona is empty for any domain with boundary contained in the graph of a $C^{1+\epsilon}$-function. This result is an extension of an earlier result of Jones and 
Garnett [GJ] for a Denjoy domain (i.e.,
a domain with boundary contained in $\Re$). 

The corona problem can be reformulated as follows, see, e.g., [Ga]:

A collection $f_{1},\dots, f_{n}$ of functions from $H^{\infty}(X)$ satisfies
the {\em corona condition} if
\begin{equation}\label{e1}
\max_{1\leq j\leq n}|f_{j}(x)|\geq\delta>0\ \ \ {\rm for\ all}\ \ \
x\in X.
\end{equation}
The corona problem being solvable (i.e., the corona is empty)
means that the Bezout equation
\begin{equation}\label{e2}
f_{1}g_{1}+\cdots+f_{n}g_{n}\equiv 1
\end{equation}
has a solution $g_{1},\dots, g_{n}\in H^{\infty}(X)$ for any 
$f_{1},\dots, f_{n}$ satisfying the corona condition. We refer to
$\max_{1\leq j\leq n}||g_{j}||_{\infty}$ as a ``bound on the corona
solutions``. (Here $||\cdot||_{\infty}$ is the norm on $H^{\infty}(X)$.)

In [Br4, Theorem 1.1] using an $L^{2}$ cohomology technique we proved 
\begin{Th}\label{te1}
Let $r:X\to Y$ be a connected unbranched covering of a Caratheodory  hyperbolic
Riemann surface of finite type $Y$ (i.e., the fundamental group of $Y$ is finitely generated). Then $X$ is Caratheodory hyperbolic and
for any $f_{1},\dots, f_{n}\in H^{\infty}(X)$ satisfying (\ref{e1}) 
there are solutions $g_{1},\dots, g_{n}\in H^{\infty}(X)$ of
(\ref{e2}) with the bound $\max_{1\leq j\leq n}||g_{j}||_{\infty}\leq 
C(Y,n,\max_{1\leq j\leq n}||f_{j}||_{\infty},\delta)$.
\end{Th}

This result extends the class of Riemann surfaces for which the corona theorem is valid (see also [Br1]). On the other hand, from the results of L\'{a}russon [L] (sharpened in [Br3]) one obtains that the assumption of the Caratheodory hyperbolicity of $Y$ cannot be removed. Specifically, for any integer $n\geq 2$ there are a compact
Riemann surface $S_{n}$ and its regular covering $r_{n}:\widetilde S_{n}\to
S_{n}$ such that
\begin{itemize}
\item[(a)]
$\widetilde S_{n}$ is a complex submanifold of an open Euclidean ball 
$\Bo_{n}\subset\Co^{n}$;
\item[(b)]
the embedding $i:\widetilde S_{n}\hookrightarrow\Bo_{n}$ induces an isometry
$i^{*}:H^{\infty}(\Bo_{n})\to H^{\infty}(\widetilde S_{n})$.
\end{itemize}
In particular, the maximal ideal spaces of 
$H^{\infty}(\widetilde S_{n})$ and $H^{\infty}(\Bo_{n})$ coincide.

The main result of our paper is the following noncommutative analog of the above theorem:
\begin{Th}\label{te2}
Let $r: X\to Y$ satisfy the assumptions of Theorem \ref{te1} and
$A=(a_{ij})$ be an $n\times k$ matrix, $k<n$, with entries in 
$H^{\infty}(X)$. Assume that the family of minors of order $k$ of $A$
satisfies the corona condition. Then there is an $n\times n$
matrix $\widetilde A=(\widetilde a_{ij})$, 
$\widetilde a_{ij}\in H^{\infty}(X)$, so that $\widetilde a_{ij}=a_{ij}$ for
$1\leq j\leq k$, $1\leq i\leq n$, and $det\ \!\widetilde A=1$.

Moreover, the corresponding norm of $\widetilde A$ is bounded by a constant depending on the norm of $A$, $\delta$ (from (\ref{e1}) for the family of minors of order $k$ of $A$), $n$ and $Y$ only.
\end{Th} 

Previously, a similar result was proved for matrices with entries in $H^{\infty}(U)$ for
domains $U\hookrightarrow X$ such that the embedding induces an injective homomorphism of the 
corresponding fundamental groups and $r(U)\subset\subset Y$,
see [Br2, Theorem 1.1].
Its proof was based on a Forelli type theorem on projections
in $H^{\infty}$ (see [Br1]) and a Grauert type theorem for
``holomorphic`` vector bundles on maximal ideal spaces (which are not usual
manifolds) of certain Banach algebras (see [Br2]) along with some ideas of Tolokonnikov [T] (see also this paper for further results and references on the extension problem for matrices with entries in different function algebras).

The remarkable class of Riemann surfaces $X$ for which a Forelli
type theorem and, hence, the corona theorem are valid was introduced by Jones and Marshall [JM]. The definition
is in terms of an interpolating property for the critical points of the Green function on $X$. 
It is an interesting open question whether the result analogous to Theorem \ref{te2} is valid for such $X$.
\sect{\hspace*{-1em}. Auxiliary Results.}
\quad {\bf 2.1.}
For a set of indices $\Lambda$ consider the family $X_{\Lambda}:=\{X_{\lambda}\}_{\lambda\in\Lambda}$ where each $X_{\lambda}$ is a connected unbranched covering of $Y$. By $r_{\lambda}:=X_{\lambda}\to Y$ we denote the corresponding projection. Considering this family as the disjoint union of sets $X_{\lambda}$ we introduce the natural complex structure on $X_{\Lambda}$. Thus $r_{\Lambda}:X_{\Lambda}\to Y$ is an unbranched covering of $Y$ where
$r_{\Lambda}|_{X_{\lambda}}:=r_{\lambda}$. 

We say that a function $f$ on $X_{\Lambda}$ belongs to $H^{\infty}(X_{\Lambda})$ if
$f|_{X_{\lambda}}\in H^{\infty}(X_{\lambda})$, $\lambda\in\Lambda$, and $\sup_{\lambda\in\Lambda}||f|_{X_{\lambda}}||_{\infty}<\infty$.
\begin{Proposition}\label{pr2.1}
The corona theorem is valid for $H^{\infty}(X_{\Lambda})$.
\end{Proposition}
{\bf Proof.} Let $f_{1},\dots, f_{n}\in H^{\infty}(X_{\Lambda})$ satisfy the corona condition (\ref{e1}). We set $f_{j\lambda}:=f_{j}|_{X_{\lambda}}$. Then each family $f_{1\lambda},\dots, f_{n\lambda}\in H^{\infty}(X_{\lambda})$ satisfies (\ref{e1}) with the same $\delta$ as for $f_{1},\dots, f_{n}$. According to Theorem \ref{te1} there are functions
$g_{1\lambda},\dots, g_{n\lambda}\in H^{\infty}(X_{\lambda})$ such that
$$
f_{1\lambda}g_{1\lambda}+\cdots+f_{n\lambda}g_{n\lambda}\equiv 1
$$
and
$$
\max_{1\leq j\leq n}||g_{j\lambda}||_{\infty}\leq C(Y,n,\max_{1\leq j\leq n}||f_{j}||_{H^{\infty}(X_{\Lambda})},\delta).
$$
Let us define $g_{1},\dots , g_{n}\in H^{\infty}(X_{\Lambda})$ by the formulas
$$
g_{j}|_{X_{\lambda}}:=g_{j\lambda}.
$$
Then $g_{1}f_{1}+\cdots+g_{n}f_{n}\equiv 1$.\ \ \ \ \ $\Box$

Let ${\cal M}_{\Lambda}$ be the maximal ideal space of the Banach algebra $H^{\infty}(X_{\Lambda})$. According to Theorem \ref{te1}, $H^{\infty}(X_{\Lambda})$ separates the points of $X_{\Lambda}$. Thus $X_{\Lambda}$ can be regarded as a subset of ${\cal M}_{\Lambda}$. Now, by Proposition \ref{pr2.1}, $X_{\Lambda}$ is dense in ${\cal M}_{\Lambda}$ in the Gelfand topology. 

We will show that Theorem \ref{te2} follows directly from
\begin{Th}\label{te2.2}
Let $A=(a_{ij})$ be an $n\times k$ matrix, $k<n$, with entries in $H^{\infty}(X_{\Lambda})$. Assume that the family of minors of order $k$ of $A$
satisfies the corona condition. Then there is an $n\times n$
matrix $\widetilde A=(\widetilde a_{ij})$, 
$\widetilde a_{ij}\in H^{\infty}(X_{\Lambda})$, so that $\widetilde a_{ij}=a_{ij}$ for
$1\leq j\leq k$, $1\leq i\leq n$, and $det\ \!\widetilde A=1$.
\end{Th}

{\bf 2.2.} We recall some constructions and results presented in [Br2]. 

According to a construction of [Br2, section 2] the covering $r_{\Lambda}:X_{\Lambda}\to Y$ can be considered as a fibre bundle over $Y$ with a discrete fibre $F_{\Lambda}$,
where $F_{\Lambda}$ is the disjoint union of the fibres $F_{\lambda}$ of the coverings $r_{\lambda}:X_{\lambda}\to Y$,
$\lambda\in\Lambda$. Let $l^{\infty}(F_{\Lambda})$ be the Banach algebra of bounded complex-valued functions $f$ on the discrete set $F_{\Lambda}$ with pointwise multiplication and norm $||f||=\sup_{x\in F_{\Lambda}}|f(x)|$. Let $\beta F_{\Lambda}$ be the {\em Stone-\v{C}ech compactification} of $F_{\Lambda}$, i.e., the maximal ideal space of $l^{\infty}(F_{\Lambda})$ equipped with the Gelfand topology.
Then $F_{\Lambda}$ is naturally embedded into $\beta F_{\Lambda}$ as an open dense subset, and the topology on $F_{\Lambda}$ induced by this embedding coincides with the original one, i.e., is discrete. Every function $f\in l^{\infty}(F_{\Lambda})$ has a unique extension $\hat f\in C(\beta F_{\Lambda})$. Further, any homeomorphism $\phi: F_{\Lambda}\to F_{\Lambda}$ determines an isometric isomorphism of Banach algebras $\phi^{*}:l^{\infty}(F_{\Lambda})\to l^{\infty}(F_{\Lambda})$. Therefore $\phi$ can be extended to a homeomorphism $\hat\phi:\beta F_{\Lambda}\to\beta F_{\Lambda}$.
From here, taking closures in $\beta F_{\lambda}$ of fibres of the bundle $r_{\Lambda}:X_{\Lambda}\to Y$, we obtain a fibre bundle $\hat r_{\Lambda}: E(Y,\beta F_{\Lambda})\to Y$ with fibre $\beta F_{\Lambda}$ so that $X_{\Lambda}$ is an open dense subset of $E(Y,\beta F_{\Lambda})$ (in fact, an open subbundle of $E(Y,\beta F_{\Lambda})$) and $\hat r_{\Lambda}|_{X_{\Lambda}}=r_{\Lambda}$. Moreover, it was proved in [Br2, Proposition 2.1] that 
\begin{itemize}
\item[(1)]
{\em for every $h\in H^{\infty}(X_{\Lambda})$ there is a unique $\hat h\in C(E(Y,\beta F_{\Lambda}))$ such that $\hat h|_{X_{\Lambda}}=h$.}
\end{itemize}
Also, it was proved in [Br4, Theorem 1.5] that for every $x\in Y$ and every $\lambda\in\Lambda$ the sequence $r_{\lambda}^{-1}(x)\subset X_{\lambda}$ is interpolating for $H^{\infty}(X_{\lambda})$ with the constant of interpolation bounded by a number depending on $x$ and $Y$ only. This immediately implies that 
\begin{itemize}
\item[(2)]
{\em for each $f\in l^{\infty}(r_{\Lambda}^{-1}(x))$ there is a function $\tilde f\in H^{\infty}(X_{\Lambda})$ such that $\tilde f|_{r_{\Lambda}^{-1}(x)}=f$.}
\end{itemize}
 In particular, the continuous extension of the algebra $H^{\infty}(X_{\Lambda})$ to $E(Y,\beta F_{\Lambda})$ separates the points on $E(Y,\beta F_{\Lambda})$. Thus $E(Y,\beta F_{\Lambda})$ can be regarded as a dense subset of ${\cal M}_{\Lambda}$.
 
Let $(U_{i})_{i\in I}$ be a countable cover of $Y$ by compact subsets
$U_{i}\subset Y$ homeomorphic to a closed ball in $\Re^{2}$. Then by our construction $\widehat U_{i}:=\hat r_{\Lambda}^{-1}(U_{i})$ is homeomorphic to $U_{i}\times\beta F_{\Lambda}$. So, $E(Y,\beta F_{\Lambda})$ is a countable union of compact subsets $\widehat U_{i}$. Since the covering dimension $dim\ \!\widehat U_{i}$ of $\widehat U_{i}$ is $2$, $i\in I$, this implies 
(cf. [Br2, Proposition 4.1])
\begin{itemize}
\item[(3)]
$$
dim\ \!E(Y,\beta F_{\Lambda})=2.
$$
\end{itemize}
Taking now an open countable cover of $Y$ by relatively compact subsets homeomorphic to an open ball in $\Re^{2}$ and the corresponding open cover of $E(Y,\beta F_{\Lambda})$ by their preimages with respect to $\hat r_{\Lambda}$ we get
\begin{itemize}
\item[(4)]
{\em $E(Y,\beta F_{\Lambda})$ is an open dense subset of ${\cal M}_{\Lambda}$, and the restriction of the Gelfand topology on ${\cal M}_{\Lambda}$ to $E(Y,\beta F_{\Lambda})$ coincides with the topology of $E(Y,\beta F_{\Lambda})$.}
\end{itemize}

{\bf 2.3.} Since $Y$ is a Riemann surface of finite type,
the theorem of Stout [St, Theorem 8.1] implies that there exist a compact
Riemann surface $R$ and a holomorphic embedding $\phi:Y\to R$ such that
$R\setminus\phi(Y)$ consists of finitely many closed disks with analytic
boundaries together with finitely many isolated points. Since $Y$ is
Caratheodory hyperbolic, the set of the disks in $R\setminus\phi(Y)$ is not
empty. Also, without loss of generality we may and will assume that
the set of isolated points in $R\setminus\phi(Y)$ is not empty, as well.
(For otherwise, $\phi(Y)$ is a bordered Riemann surface and the required
result follow from [Br2, Theorem 1.1].) We will naturally identify
$Y$ with $\phi(Y)$. Also, we set 
\begin{equation}\label{eq2.1}
R\setminus Y:=\left(\bigsqcup_{1\leq i\leq k}\overline{D}_{i}\right)\bigcup
\left(\bigcup_{1\leq j\leq l}\{x_{j}\}\right)\ \ \ {\rm and}\ \ \
Z:=Y\bigcup\left(\bigcup_{1\leq j\leq l}\{x_{j}\}\right)
\end{equation}
where each $D_{i}$ is biholomorphic to the open unit disk $\Di\in\Co$
and these biholomorphisms are extended to diffeomorphisms of the closures
$\overline{D}_{i}\to\overline\Di$. Then $Z\subset R$ is a bordered Riemann surface with a nonempty boundary. In particular, there is a bordered Riemann surface $Z'$ containing $\overline{Z}$ such that $\overline{Z}$ is a deformation retract of $Z'$. We set
\begin{equation}\label{eq2.2}
Y':=Z'\setminus\{x_{1},\dots,x_{l}\}.
\end{equation}
Then $Y\subset Y'$ and $\pi_{1}(Y)\cong\pi_{1}(Y')$ (here
$\pi_{1}(M)$ stands for the fundamental group of $M$). This implies that for each $\lambda\in\Lambda$ there is a connected covering
$X'_{\lambda}$ of $Y'$ such that $X_{\lambda}$ is an open connected subset of $X_{\lambda}'$. Without loss of generality we denote the covering projection $X'_{\lambda}\to Y'$ by the same symbol $r_{\lambda}$ (as for $X_{\lambda}$).  Now, we define
$X_{\Lambda}':=\{X_{\lambda}'\}_{\lambda\in\Lambda}$ so that
$X_{\Lambda}$ is an open subset of $X_{\Lambda}'$ and 
$r_{\Lambda}:X_{\Lambda}'\to Y'$, $r_{\Lambda}|_{X_{\lambda}'}:=r_{\lambda}$. 

Further, similarly to the constructions of section 2.2 we determine the bundle
$\hat r_{\Lambda}: E(Y',\beta F_{\Lambda})\to Y'$ so that
$E(Y,\beta F_{\Lambda})$ is an open subbundle of $E(Y',\beta F_{\Lambda})$. Then $X_{\Lambda}'$ and $E(Y',\beta F_{\Lambda})$ possess the properties similar to (1)-(3) for 
$X_{\Lambda}$ and $E(Y,\beta F_{\Lambda})$.

Let $cl(Y)$ denote the closure of $Y$ in $Y'$. We set
$$
E(cl(Y),\beta F_{\Lambda}):=\hat r_{\Lambda}^{-1}(cl(Y)).
$$
Then we have
\begin{itemize}
\item[(5)]
{\em $dim\ \!E(cl(Y),\beta F_{\Lambda})=2$ and
$E(Y,\beta F_{\Lambda})\subset E(cl(Y),\beta F_{\Lambda})$ is an open dense subset.}
\end{itemize}

{\bf 2.4.} By $H^{\infty}(E(Y,\beta F_{\Lambda}))$ we denote the extension of $H^{\infty}(X_{\Lambda})$ to $E(Y,\beta F_{\Lambda})$ described in section 2.2. We will use also the
algebra $H^{\infty}(E(Y',\beta F_{\Lambda}))$ determined similarly (i.e., with $Y$ and $X_{\Lambda}$ substituted for $Y'$ and $X_{\Lambda}'$).

Next, let us consider Banach subalgebras ${\cal A}_{1}$, ${\cal A}_{2}$ of $H^{\infty}(E(Y,\beta F_{\Lambda}))$ defined as follows.
\begin{equation}\label{eq2.3}
{\cal A}_{1}:=\{\hat r_{\Lambda}^{*}f\in H^{\infty}(E(Y,\beta F_{\Lambda}))\ :\
f\in H^{\infty}(Z')\}. 
\end{equation}
(Here $\hat r_{\Lambda}^{*}f$ is the pullback of $f$ with respect to $\hat r_{\Lambda}$.) 

To define ${\cal A}_{2}$ we choose a function $\phi\in H^{\infty}(Z')$ with the set of zeros $\{x_{1},\dots, x_{l}\}$ so that each $x_{j}$ is a zero of order $1$ of $\phi$. (Since
$Z'\subset\subset R$ is a bordered Riemann surface with a nonempty boundary, such a $\phi$ exists due to [Br2, Corollary 1.8].)
Then ${\cal A}_{2}$ is the uniform closure of the algebra of functions $f\in H^{\infty}(E(Y,\beta F_{\Lambda}))$ of the form
\begin{equation}\label{eq2.4}
f:=g+([\hat r_{\Lambda}^{*}\phi]\cdot h)|_{E(Y,\beta F_{\Lambda})},\ g\in {\cal A}_{1},\ h\in H^{\infty}(E(Y',\beta F_{\Lambda})).
\end{equation}
By the definition ${\cal A}_{2}$ separates the points of $E(Y,\beta F_{\Lambda})$ (because $H^{\infty}(E(Y',\beta F_{\Lambda}))$ separates the points of $E(Y',\beta F_{\Lambda})$
and $\hat r_{\Lambda}^{*}\phi$ is nonzero on the fibres of $\hat r_{\Lambda}$).

Clearly, we have embeddings
\begin{equation}\label{eq2.5}
{\cal A}_{1}\stackrel{i_{1}}{\longrightarrow}{\cal A}_{2}\stackrel{i_{2}}{\longrightarrow}H^{\infty}(E(Y,\beta F_{\Lambda})).
\end{equation}
The transpose maps to these embeddings determine continuous surjective maps
\begin{equation}\label{eq2.6}
{\cal M}_{\Lambda}\stackrel{i_{2}^{*}}{\longrightarrow}M_{2}\stackrel{i_{1}^{*}}{\longrightarrow}M_{1}
\end{equation}
where $M_{2}$ is the closure in the Gelfand topology of the image of $E(Y,\beta F_{\Lambda})$ in the maximal ideal space of
${\cal A}_{2}$, and $M_{1}$ is the closure in the Gelfand topology of the image of $E(Y,\beta F_{\Lambda})$ in the maximal ideal space of ${\cal A}_{1}$. (Here we used that the closure in the Gelfand topology of $E(Y,\beta F_{\Lambda})\subset {\cal M}_{\Lambda}$ is ${\cal M}_{\Lambda}$, see Proposition \ref{pr2.1}.) 

By the definition,
$M_{1}=\overline{Z}$ and $E(cl(Y),\beta F_{\Lambda})\subset M_{2}$ (see section 2.3). Moreover, the restriction of $i_{2}^{*}$ to $E(Y,\beta F_{\Lambda})$ is the identity map and the restriction of $i_{1}^{*}$ to $E(cl(Y),\beta F_{\Lambda})$ can be naturally identified with $\hat r_{\Lambda}$ so that
$(i_{1}^{*})^{-1}(cl(Y))=\hat r_{\Lambda}^{-1}(cl(Y))=E(cl(Y),\beta F_{\Lambda})$.
Now, we prove
\begin{Lm}\label{le2.3}
For each $x_{j}\in M_{1}$ the compact set $(i_{1}^{*})^{-1}(x_{j})$ consists of a single point (which we naturally identify with $x_{j}$), $1\leq j\leq l$.
\end{Lm}
{\bf Proof.} Let $\{\xi_{1,\alpha}\}, \{\xi_{2,\alpha}\}\subset E(Y,\beta F_{\Lambda})$ be nets converging to points $\xi_{1}, \xi_{2}\in (i_{1}^{*})^{-1}(x_{j})$. Then for $f$ from (\ref{eq2.4}) and $i=1,2$ we have
\begin{equation}\label{eq2.7}
f(\xi_{i})=\lim_{\alpha}f(\xi_{i,\alpha})=
\lim_{\alpha}(g(\xi_{i,\alpha})+(\hat r_{\Lambda}^{*}\phi)(\xi_{i,\alpha})\cdot h(\xi_{i,\alpha})):=g(x_{j}).
\end{equation}
(We used here that the nets $\{i_{1}(\xi_{1,\alpha})\}, \{i_{1}(\xi_{2,\alpha})\}\subset M_{1}$ converge to $x_{j}$.)

This implies that $\xi_{1}=\xi_{2}$.\ \ \ \ \ $\Box$
\begin{C}\label{cor2.4}
$$
dim\ \! M_{2}=2.
$$
\end{C}
{\bf Proof.} According to Lemma \ref{le2.3} and property (5) of section 2.3, $M_{2}$ is the disjoint union of zero-dimensional sets $\{x_{j}\}$, $1\leq j\leq l$, and the two-dimensional set
$E(cl(Y),\beta F_{\Lambda})$. Hence $dim\ \! M_{2}=2$, see, e.g., [N, Chapter 2, Theorem 9-11].\ \ \ \ \
$\Box$\\

{\bf 2.5.} We fix coordinate neighbourhoods $N_{j}\subset\subset Z$ (see (\ref{eq2.1})) biholomorphic to $\Di$ of points $x_{j}$, $1\leq j\leq l$, and a bordered Riemann surface $S\subset Y$ such that $N_{i}\cap N_{j}=\emptyset$ for $i\neq j$, each $Y\cap N_{j}$ does not contain $x_{j}$ and is biholomorphic to an annulus and
${\cal U}:=S\cup(\cup_{1\leq j\leq l}N_{j}^{*})$ is an open cover of $Y$. Here $N_{j}^{*}:=N_{j}\setminus\{x_{j}\}$ is biholomorphic to $\Di^{*}:=\Di\setminus\{0\}$. We set
\begin{equation}\label{eq2.8}
N_{j\Lambda}^{*}:=r_{\Lambda}^{-1}(N_{j}^{*}),\ \ \ 1\leq j\leq l,\ \ \ {\rm and}\ \ \ S_{\Lambda}:=r_{\Lambda}^{-1}(S).
\end{equation}

Let $V\subset X_{\Lambda}$ be either one of $N_{j\Lambda}^{*}$ or $S_{\Lambda}$.
By $H^{\infty}(V)$ we denote the Banach algebra of bounded holomorphic functions on $V$ defined as in section 2.1 for $X_{\Lambda}$. Further, we set
\begin{equation}\label{eq2.9}
\hat N_{j}:=(i_{2}^{*}\circ i_{1}^{*})^{-1}(N_{j}),\ \ \ 1\leq j\leq l,\ \ \ \hat S:=(i_{2}^{*}\circ i_{1}^{*})^{-1}(S\cup\partial Z)
\end{equation}
Here $\partial Z$ is the boundary of the bordered Riemann surface $Z$ that can be regarded as the ``outer boundary`` of $S$.

By the definition $\hat N_{j}$, $1\leq j\leq l$, and $\hat S$ are open subsets of ${\cal M}_{\Lambda}$ forming a cover of this space.
The main fact used in the proof of Theorem \ref{te1} is
\begin{Proposition}\label{pr2.5}
Assume that $f\in H^{\infty}(V)$ where $V$ is either one of $N_{j\Lambda}^{*}$ or $S_{\Lambda}$. Then $f$ admits a continuous extension $\hat f$ to $\hat V$ where $\hat V$ stands for the corresponding $\hat N_{j}$ or $\hat S$.
\end{Proposition}
{\bf Proof.} First, we will prove the result for $V=N_{j\Lambda}^{*}$. 

Let $\rho_{j}$ be a $C^{\infty}$-function on $R$ equal to $1$ in a neighbourhood of $x_{j}$ with $supp(\rho)\subset\subset N_{j}$. 
We set
\begin{equation}\label{eq2.10}
f_{1}:=(r_{\Lambda}^{*}\rho_{j})\cdot f.
\end{equation}
Then $f_{1}$ can be considered as a $C^{\infty}$-function on $X_{\Lambda}'$ (defined in section 2.3). Further, we introduce a $(0,1)$-form on $X_{\Lambda}'$ by the formula
\begin{equation}\label{eq2.11}
\omega:=\frac{\overline\partial f_{1}}{\rho_{\Lambda}^{*}\phi}.
\end{equation}
The definition is correct because $\overline\partial f_{1}$ equals $0$ on $\rho_{\Lambda}^{-1}(O)$ for some neighbourhood $O$ of $x_{j}$ and on $X_{\Lambda}'\setminus N_{j\Lambda}^{*}$, and $\rho_{\Lambda}^{*}\phi\neq 0$ on $N_{j\Lambda}^{*}$. Thus 
$\omega$ is a $\overline\partial$-closed $1$-form on $X_{\Lambda}'$. Consider the form 
$$
\omega_{\lambda}:=\omega|_{X_{\lambda}'}\ \ \ {\rm on}\ \ \ X_{\lambda}'.
$$

Let us assume that $Z'$ is equipped with a hermitian metric $h_{Z'}$ with
the associated $(1,1)$-form $\omega_{Z'}$. Then we equip  $X_{\lambda}'$ with the
hermitian metric $h_{X_{\lambda}'}$ induced by the pullback 
$r_{\lambda}^{*}\omega_{Z'}$ of $\omega_{Z'}$ to $X_{\lambda}'$. 
Now, if $\eta$ is a smooth $(0,1)$-form on $X_{\lambda}'$, by
$|\eta|_{z}$, $z\in X_{\lambda}'$, we denote the norm of $\eta$ at $z$
defined by the hermitian metric $h_{X_{\lambda}'}^{*}$ on the fibres of the cotangent
bundle $T^{*}X_{\lambda}'$ on $X_{\lambda}'$. 

Next, since $f\in  H^{\infty}(N_{j\Lambda}^{*})$ and 
$r_{\Lambda}(supp(\omega))=:K\subset\subset Y'$, see (\ref{eq2.2}), one easily obtains from
(\ref{eq2.11}) that
\begin{equation}\label{eq2.12}
||\omega||:=\sup_{\lambda\in\Lambda}\left\{\sup_{z\in X_{\lambda}}|\omega_{\lambda}|_{z}\right\}<\infty
\end{equation}
From here by [Br4, Theorem 1.6] we obtain that
the equation
$\overline\partial g_{\lambda}=\omega_{\lambda}$
has a smooth bounded solution $g_{\lambda}$ on $X_{\lambda}'$ such that
\begin{equation}\label{eq2.13}
||g_{\lambda}||_{L^{\infty}}:=\sup_{z\in X_{\lambda}'}|g_{\lambda}(z)|\leq C||\omega||
\end{equation}
with $C$ depending on $K$, $Z'$ and $h_{Z'}$ only.

We define bounded functions $g$ and $f_{2}$ on $X_{\Lambda}'$ by the formulas
\begin{equation}\label{eq2.14}
g|_{X_{\lambda}'}:=g_{\lambda},\ \ \ \lambda\in\Lambda,\ \ \
f_{2}:=(\rho_{\Lambda}^{*}\phi)\cdot g.
\end{equation}
Then we have 
\begin{equation}\label{eq2.15}
(a)\ \ \overline\partial f_{2}=\overline\partial f_{1}\ \ \ {\rm on}\ \ \ X_{\Lambda}'\ \ \ {\rm and}\ \ \
(b)\ \ \lim_{\alpha}f_{2}(\xi_{\alpha})=0
\end{equation}
for each net $\{\xi_{\alpha}\}\subset X_{\Lambda}'$ such that
$\{r_{\Lambda}(\xi_{\alpha})\}\subset Y'$ is a net converging to
any $x_{s}$, $1\leq s\leq l$. In particular,
\begin{equation}\label{eq2.16}
f_{3}:=f_{1}-f_{2}\in H^{\infty}(X_{\Lambda}').
\end{equation}
Thus $f_{3}$ admits a continuous extension $\hat f_{3}$ to ${\cal M}_{\Lambda}$.

Let us prove now that
\begin{itemize}
\item[(*)]
{\em
$f_{2}$ admits a continuous extension $\hat f_{2}$ to ${\cal M}_{\Lambda}$.}
\end{itemize}

Indeed, by the definition of $f$ and $r_{\Lambda}^{*}\rho_{j}$, the function $f_{1}$ defined by (\ref{eq2.10}) has a continuous extension to $E(Y',\beta F_{\Lambda})$, see [Br2, Proposition 2.1]. Thus, $f_{2}:=f_{1}-f_{3}$ admits a continuous extension to $E(Y',\beta F_{\Lambda})$, as well. (We denote this extension also by $f_{2}$.) Now if $\{\xi_{\alpha}\}\subset E(Y',\beta F_{\Lambda})$ is a net converging to a point $\xi\in M_{2}$ (see (\ref{eq2.6})) such that $i_{1}^{*}(\xi)=x_{s}$ for some $1\leq s\leq l$, then from (\ref{eq2.15}) (b) we get
$$
\lim_{\alpha}f_{2}(\xi_{\alpha})=0.
$$
Since $(i_{1}^{*})^{-1}(x_{s})=x_{s}$, the latter implies that the function $\tilde f_{2}$ equals $0$ at each $x_{s}$ and $f_{2}$ on $E(cl(Y),\beta F_{\Lambda})$ is continuous on $M_{2}$. Therefore the function $\hat f_{2}:=i_{2}^{*}\tilde f_{2}$ is continuous on ${\cal M}_{\Lambda}$. Since the restriction of $i_{2}^{*}$ to $E(Y,\beta F_{\Lambda})$ is the identity map, $\hat f_{2}$ is a continuous extension of $f_{2}$.
This proves (*).

From (\ref{eq2.16}) and (*) we obtain that $f_{1}$ admits a continuous extension $\hat f_{1}$ to ${\cal M}_{\Lambda}$.
Now, $\hat N_{j}$ from (\ref{eq2.9}) is the union of 
$\hat r_{\Lambda}^{-1}(N_{j}^{*})\subset E(Y,\beta F_{\Lambda})$ and $(i_{2}^{*}\circ i_{1}^{*})^{-1}(O)$ where
$O\subset\subset N_{j}$ is a neighbourhood of $x_{j}$ such that
$\rho_{j}\equiv 1$ on $O$. The function $f$ admits a continuous extension $\tilde f$ on $\hat r_{\Lambda}^{-1}(N_{j}^{*})$, see
[Br2, Proposition 2.1], and $f=f_{1}$ on $r_{\Lambda}^{-1}(O)$.
Thus the function $\hat f$ defined by
$$
\hat f:=\hat f_{1}\ \ \ {\rm on}\ \ \ (i_{2}^{*}\circ i_{1}^{*})^{-1}(O)\ \ \ {\rm and}\ \ \ \hat f:=\tilde f\ \ \
{\rm on}\ \ \ \hat r_{\Lambda}^{-1}(N_{j}^{*})
$$
is the required continuous extension of $f$ to $\hat N_{j}$.

Finally, in the case $V=S_{\Lambda}$ we choose a $C^{\infty}$-function $\rho$ on $R$ equals $0$ on $Y\setminus S$ and $1$ on $R\setminus Z$ with $supp(d\rho)\subset\subset S$.
Then repeating the above arguments with $\rho_{j}$ substituted for $\rho$ we obtain the proof of the proposition in this case.
We leave the details to the readers.\ \ \ \ \ $\Box$
\sect{\hspace*{-1em}. Proof of Theorem \ref{te2}.}
\quad {\bf 3.1. Proof of Theorem \ref{te2.2}.} Let $A=(a_{ij})$ be an $n\times k$ matrix, $k<n$, with entries in $H^{\infty}(X_{\Lambda})$. Assume that the family of minors of order $k$ of $A$ satisfies the corona condition (\ref{e1}).
Due to the corona theorem for $H^{\infty}(X_{\Lambda})$, see Proposition \ref{pr2.1}, we can extend $A$ continuously to ${\cal M}_{\Lambda}$ such that the family of minors of order $k$ of the extended matrix $\hat A=(\hat a_{ij})$ satisfies (\ref{e1}) on ${\cal M}_{\Lambda}$ with the same $\delta$ as for $A$. Next, according to [L, Theorem 3], to prove the theorem it suffices to find  an $n\times n$ matrix $B=(b_{ij})$, $b_{ij}\in C({\cal M}_{\Lambda})$, so that $b_{ij}=\hat a_{ij}$ for
$1\leq j\leq k$, $1\leq i\leq n$, and $det\ \!B=1$.

Note that the matrix $\hat A$ determines a trivial subbundle
$\xi$ of complex rank $k$ in the trivial vector bundle $\theta^{n}:={\cal M}_{\Lambda}\times\Co^{n}$ on ${\cal M}_{\Lambda}$. Let $\eta$ be an additional to $\xi$ subbundle of $\theta^{n}$, i.e., $\xi\oplus\eta=\theta^{n}$. We will prove that $\eta$ is topologically trivial. Then a  trivialization $s_{1},s_{2},\dots, s_{n-k}\in C({\cal M}_{\Lambda},\eta)$ (given by global continuous sections of $\eta$) will determine the required continuous extension $B$ of the matrix $\hat A$.

Let us prove first that $\hat A$ can be extended to an invertible matrix on each $\hat N_{j}$ and $\hat S$, see (\ref{eq2.9}).
\begin{Lm}\label{le3.1}
Let $\hat V$ be either one of $\hat N_{j}$ or $\hat S$. Then
for $\hat A|_{\hat V}$ there is an $n\times n$ matrix $B_{\hat V}=(b_{ij;\hat V})$, $b_{ij;\hat V}\in C(\hat V)$, so that $b_{ij;\hat V}=\hat a_{ij}|_{\hat V}$ for
$1\leq j\leq k$, $1\leq i\leq n$, and $det\ \!B_{\hat V}=1$.
Moreover, $B_{\hat V}|_{V}$ has entries in $H^{\infty}(V)$ where
$V:=\hat V\cap X_{\Lambda}$.
\end{Lm}
{\bf Proof.} First assume that $\hat V=\hat N_{j}$ so that
$V=N_{j\Lambda}^{*}$ is an unbranched covering of $N_{j}^{*}\cong\Di^{*}$, see (\ref{eq2.8}). Then by the definition $N_{j\Lambda}^{*}=\{N_{j\lambda}^{*}\}_{\lambda\in\Lambda}$ where each $N_{j\lambda}^{*}:=N_{j\Lambda}^{*}\cap X_{\lambda}$ is an unbranched covering of $N_{j}^{*}$ consisting of at most countably many connected components. Thus each $N_{j\lambda}^{*}$ is biholomorphic to $\sqcup_{k\in K_{\lambda}}W_{j\lambda; k}$, $K_{\lambda}\subset\N$, where each $W_{j\lambda;k}$ is either $\Di$ or $\Di^{*}$. Now, $A|_{W_{j\lambda;k}}$ satisfies
conditions of Theorem \ref{te1} with the same $\delta$ as for $A$. According to the main result of Tolokonnikov [T] for $H^{\infty}$-matrices on $\Di$, there is a matrix $B_{j\lambda;k}$ with entries in $H^{\infty}(W_{j\lambda;k})$ which extends $A|_{W_{j\lambda;k}}$ in the sense of Theorem \ref{te1} and such that $det\ \! B_{j\lambda;k}=1$ and
\begin{equation}\label{eq3.1}
\sup_{j,\lambda, k}||B_{j\lambda;k}||\leq C
\end{equation}
where $C$ depends on the norm of $A$ on $X_{\Lambda}$, $\delta$ and $n$.
(Here for a matrix $C=(c_{ij})$ with entries in $H^{\infty}(O)$ we set $||C||:=\max_{i,j}||c_{ij}||_{H^{\infty}(O)}$.)
In particular, (\ref{eq3.1}) implies that the matrix $B_{j\Lambda}$ on $N_{j\Lambda}^{*}$ defined by
$$
B_{j\Lambda}|_{W_{j\lambda;k}}=B_{j\lambda;k},\ \ \ 1\leq j\leq l,\ k\in K_{\lambda},\ \lambda\in\Lambda ,
$$
has entries in $H^{\infty}(N_{j\Lambda}^{*})$, extends $A|_{N_{j\Lambda}^{*}}$ and $det\ \! B_{j\Lambda}=1$. According to Proposition \ref{pr2.5}, $B_{j\Lambda}$ is extended to a continuous matrix $B_{\hat N_{j}}$ on $\hat N_{j}$. This matrix extends $\hat A|_{\hat N_{j}}$ and satisfies the required conditions of the lemma.

Consider now the case $\hat V=\hat S$ so that $V=S_{\Lambda}$ is an unbranched covering of $S$, see (\ref{eq2.8}). In this case we apply similar to the above arguments where instead of the result of [T] we use [Br2, Theorem 1.1] applied to the coverings $S_{\lambda}:=S_{\Lambda}\cap X_{\lambda}$ of a bordered Riemann surface $S$. Then we obtain a matrix $B_{\Lambda}$ on $S_{\Lambda}$ with entries in $H^{\infty}(S_{\Lambda})$ which extends $A|_{S_{\Lambda}}$ and such that $det\ \! B_{\Lambda}=1$. Applying again Proposition \ref{pr2.5}, we extend $B_{\Lambda}$ continuously to $\hat S$ so that the extended matrix $B_{\hat S}$ satisfies the required conditions of the lemma.\ \ \ \ \ $\Box$

Let $\xi_{q}$ be the quotient bundle of $\theta^{n}$ with respect to the subbundle $\xi$. By the definition $\xi_{q}$ is isomorphic (in the category of continuous bundles on ${\cal M}_{\Lambda}$) to $\eta$. Thus it suffices to prove that $\xi_{q}$ is topologically trivial. 

Now, Lemma \ref{le3.1} implies straightforwardly that $\xi_{q}|_{\hat V}$ is topologically trivial for $\hat V$ being either one of $\hat N_{j}$ or $\hat S$. In particular, $\xi_{q}$ is defined by a $1$-cocycle defined on the open cover
$\{\hat N_{1},\dots ,\hat N_{l}, \hat S\}$ of ${\cal M}_{\Lambda}$ (see, e.g., [H] for the general theory of vector bundles). Since by the definition $\hat N_{i}\cap\hat N_{j}=\emptyset$ for $i\neq j$, this cocycle consists of continuous matrix-functions
$$
C_{i}\in C(\hat N_{i}\cap \hat S,GL_{n-k}(\Co)),\ \ \
1\leq i\leq l.
$$

Set now $\tilde N_{j}:=(i_{1}^{*})^{-1}(N_{j})$, $\tilde S:=(i_{1}^{*})^{-1}(S)$, see (\ref{eq2.6}). Then $\{\tilde N_{1},\dots, \tilde N_{l},\tilde S\}$ is an open cover of $M_{2}$. Moreover, the map $i_{2}: {\cal M}_{\Lambda}\to M_{2}$ is identity on each $\hat N_{i}\cap\hat S$, see section 2.4. Therefore each $C_{i}$ can be regarded as a matrix-function on $\tilde N_{i}\cap\tilde S$. In particular, these functions determine a complex vector bundle $\tilde\xi_{q}$ of rank $n-k$ on $M_{2}$ so that 
\begin{equation}\label{eq3.2}
i_{2}^{*}\tilde\xi_{q}=\xi_{q}.
\end{equation}

Since $dim\ \!M_{2}=2$, see Corollary \ref{cor2.4}, the bundle $\tilde\xi_{q}$ is isomorphic to $\theta_{M_{2}}^{n-k-1}\oplus\theta$ where $\theta_{M_{2}}^{n-k-1}:=M_{2}\times\Co^{n-k-1}$ is the trivial bundle and $\theta$ is a vector bundle of complex rank $1$, see, e.g., [Br5, Lemma 2.8]. This and (\ref{eq3.2}) imply that
$\xi_{q}\cong\eta$ is isomorphic to $\theta^{n-k-1}\oplus i_{2}^{*}\theta$ where $\theta^{n-k-1}={\cal M}_{\Lambda}\times\Co^{n-k-1}$ is the trivial bundle. Now, for the first Chern classes (which are additive with respect to the operation of the direct sum of bundles) we have the following identity
\begin{equation}\label{eq3.3}
0=c_{1}(\theta^{n})=c_{1}(\xi)+c_{1}(\eta)=c_{1}(\theta^{n-k-1}\oplus i_{2}^{*}\theta)=c_{1}(i_{2}^{*}\theta).
\end{equation}
We used here that Chern classes of trivial bundles are zeros.

Equality (\ref{eq3.3}) shows that the first Chern class of the complex rank $1$ vector bundle $i_{2}^{*}\theta$ is zero. Thus this bundle is topologically trivial (see, e.g., [H]). Combining this fact with the above isomorphism for $\eta$ we get
$\eta\cong\theta^{n-k}:={\cal M}_{\Lambda}\times\Co^{n-k}$.

This completes the proof of Theorem \ref{te2.2}.\ \ \ \ \ $\Box$
\\

{\bf 3.2. Proof of Theorem \ref{te1}.} Let us define  $\Lambda_{Y;M,\delta}$ as the set of all possible couples $(A,X)$ where
$X$ is a connected covering of $Y$ and $A$ is an $n\times k$ matrix on $X$ satisfying conditions of Theorem \ref{te1} with a fixed $\delta$ in the corona condition (\ref{e1}) for the family of minors of order $k$ and such that $||A||\leq M$. 
For $\Lambda:=\Lambda_{Y;M,\delta}$ we consider the $n\times k$ matrix ${\cal A}$ with entries in $H^{\infty}(X_{\Lambda})$ defined as follows
\begin{equation}\label{eq3.4}
{\cal A}|_{X_{\lambda}}:=A,\ \ \ \lambda=(A,X)\in\Lambda,\ \ \ X_{\lambda}:=X .
\end{equation}
Then clearly ${\cal A}$ satisfies conditions of Theorem \ref{te2.2} on $X_{\Lambda}$. According to this theorem there is an $n\times n$ matrix $\widetilde {\cal A}$ with entries in $H^{\infty}(X_{\Lambda})$ and with $det\ \!\widetilde {\cal A}=1$ that extends ${\cal A}$.
For $\lambda=(A,X)\in\Lambda$ we set
$$
\widetilde A:=\widetilde {\cal A}|_{X}.
$$
Then $\widetilde A$ extends $A$ and $\det\ \!\widetilde A=1$, and $||\widetilde A||\leq C(||A||,\delta,M,Y)$.

The proof of Theorem \ref{te1} is complete.\ \ \ \ \ $\Box$

\end{document}